\documentclass{article}
\usepackage{fullpage}
\usepackage{amsmath}
\usepackage{amsfonts}
\usepackage{amsthm}
\usepackage{amssymb}
\usepackage{mathrsfs}
\usepackage{xr-hyper}
\usepackage{color,xcolor}
\usepackage{ifpdf}
\usepackage{psfrag}
\usepackage{graphicx,graphics,subfigure,epsfig}
\usepackage{pgf,tikz}
\usepackage{url}
\usepackage{hyperref}
\usepackage{xspace}
\usepackage{bm}
\usepackage{nicefrac}
\usepackage{bbm}

\usepackage{pgf,tikz}

\newcommand{\dpar}[2]{\dfrac{\partial #1}{\partial #2}}
\newtheorem{theorem}{Theorem}[section]
\newtheorem{remark}[theorem]{Remark}

\newcommand{\R}{\mathbb R}

\renewcommand{\P}{\mathbb P}

\newcommand{\DD}{\mathcal D}

\newcommand{\FF}{\mathcal F}

\newcommand{\bba}{\mathbf{a}}
\newcommand{\bbb}{\mathbf{b}}

\newcommand{\bbf}{\mathbf{f}}

\newcommand{\bbn}{\mathbf{n}}

\newcommand{\bbu}{\mathbf{u}}
\newcommand{\bbv}{\mathbf{v}}
\newcommand{\bbw}{\mathbf{w}}
\newcommand{\bbx}{\mathbf{x}}
\newcommand{\bby}{\mathbf{y}}
\newcommand{\bbz}{\mathbf{z}}

\newcommand{\bbJ}{\mathbf{J}}
\newcommand{\bbK}{\mathbf{K}}

\newcommand{\bbN}{\mathbf{N}}

\newcommand{\bsa}{\bm\alpha}

\newcommand*\xbar[1]{%
  \hbox{%
    \vbox{%
      \hrule height 0.5pt % The actual bar
      \kern0.4ex%         % Distance between bar and symbol
      \hbox{%
        \kern-0.05em%      % Shortening on the left side
        \ensuremath{#1}%
        \kern-0.00em%      % Shortening on the right side
      }%
    }%
  }%
}

\newcommand{\dt}{\Delta t}

\newcommand{\hbbf}{\hat{\mathbf{f}}}

\newcommand{\Id}{\text{Id }}
\begin{document}
\title{Bound preserving {P}oint-{A}verage-{M}oment {P}olynomi{A}l-interpreted ({PAMPA}) on polygonal meshes}
\author{R\'emi Abgrall$^{(\star)}$,  Yongle Liu$^{(\star)}$ and Walter Boscheri $^{(\dagger)}$ 
\\
$(\star)$: Institut f\"ur Mathematik, Universit\"at Z\"urich, Z\"urich, Switzerland\\
$(\dagger)$: LAMA, Universit\'e de Chamb\'ery, Chamb\'ery, France\\
email: \{remi.abgrall,yongle.liu\}\@math.uzh.ch, walter.boscheri@univ-smb.fr}
\maketitle

\begin{abstract}
We present a novel discretisation strategy, strongly inspired from Roe's Active Flux scheme. It can use polygonal meshes and is provably bound preserving for scalar problems and the Euler equations. Several cases demonstrates the quality of the method, and improvements with respect to previous work of the authors. This paper is a summary of \cite{BPPampa}.
\end{abstract}
%\keywords{Active flux, polygonal meshes, virtual elements, bound preserving,}
%\bodymatter
\section{Introduction}
We are interested in the numerical approximation of non linear hyperbolic conservation laws, scalar and systems with  to initial and boundary conditions.
 There is already a lot of methods that have successfully achieved their discretisation. One may mention finite difference methods, finite volume methods and their  high order extensions, namely ENO, WENO % \cite{Ami,WENO} 
  and the discontinuous Galerkin methods. %\cite{CockburnKShu}. 
 Another class of method, also very successful though less popular are the finite element methods with globally continuous approximations. % \cite{mallet,edge,Guermond,KuzminHennes,DeconinckMario}. 
  Each have their advantages and drawback. For example, finite difference methods need Cartesian meshes, or image of Cartesian meshes, so that meshing complicated geometries can be difficult, though there exists reliable strategies to solve this issue, such as meshes by block. Meshing complex geometries can be addressed by using unstructured meshes, though constructing a mesh for a complex 3D object is still a challenge. In case of unstructured meshes, method of choice can be discontinuous Galerkin methods or finite element ones\footnote{About continuous finite element method, an urban legend claim they are not locally conservative. In \cite{abgr}, we show that one can construct explicitly numerical flux-though non standard-, so that the scheme can be rewritten equivalently as a finite volume method. This is also true for any method, except the one we discuss in this paper.}. It is also possible to do mesh refinement there, though this is not an easy task. Discontinuous Galerkin methods will often need the introduction of hanging nodes, and continuous FEM will need techniques that completely change the topology of the mesh. Nothing is impossible, anything has a price.

In a series of AIAA paper and PhD thesis from the University of Michigan, P.L Roe and his students have introduced a new method, that they call "Active Flux"  \cite{AF1,AF2,AF3,RoeAF}. This can be seen as a multidimensional extension of van Leer V's method  described in \cite{VLeV}. In two dimensions, with triangular meshes, see \cite{He}, the solution is described by average values in each triangle, point values at the vertices and at the mid point of the edges. The update of the average value is more or less done as in the present paper, by simply using the conservative formulation of the problem, but the flux are interpolated in time (this explains the vocable "active"). The point values, where conservation is not needed, are evolved by the method of characteristics. The method is formally third order in space, and have excellent properties for smooth problems. In our opinion, the extension to system is more complex, but this method is the initial source of inspiration of the current paper.

This approach has also inspired several authors. In \cite{HKS}, the active flux methodology is combined with the ADER time stepping technique. In \cite{zbMATH07818681}, time stepping is again addressed, this time using bicharacteristics for systems. In \cite{zbMATH07698889}, mesh adaptation is explored via an AMR strategy. In \cite{Abgrall_camc}, we show how to use the method of lines, and show that the active flux strategy is locally conservative, though the interpretation as a "pure" finite volume method is unclear. Extension to arbitrary high order is addressed in one dimension in \cite{AB_HOAF}, and very recently in \cite{barsukow2025generalizedactivefluxmethod} for Cartesian meshes, without non linear stabilisation. Interesting structure preserving properties of this new class of scheme has been described in \cite{zbMATH07695228}, but it is possible that this is the structure of the method itself, and not the type of mesh, that helps to provide these properties, see \cite{abgrall2024virtualfiniteelementhyperbolic}, all this is still very obscure.

We have been interested in developing this method for unstructured meshes in \cite{abgrall2024activefluxtriangularmeshes}, and then \cite{abgrall2024virtualfiniteelementhyperbolic} for general polygons.  In both cases, the approximation is globally continuous, as in any active flux method. In the second reference, we are interested in polygons in the aim of facilitating mesh refinement. We use a functional approximation of the data (average values and point values) that is inspired of the Virtual Finite Element literature, so that it is possible to cut edges and faces without  destroying global continuity. The global continuity allows to minimize the number of degrees of freedom. The two combined are, in our opinion, an interesting property.

The rest of the paper is devoted to describe an improvement of \cite{abgrall2024virtualfiniteelementhyperbolic} where the non linear stabilisation for discontinuities is inspired by MOOD \cite{Mood1}. Here, following \cite{wu2023geometric}, we show how to construct bound preserving scheme on polygons, without using any polynomial approximation: we simply use the data in hand. This is an extension to multiD of \cite{abgrall2024boundpreservingpointaveragemomentpolynomialinterpreted}. 

%%%%%%%%%
\bigskip
The hyperbolic problems we consider in this paper are written in a conservative form as:
\begin{equation}\label{eq:conservative}
\dpar{\bm u}{t}+\text{ div }\bm f(\bm u)=0, \quad \bm x \in \Omega\subset \R^d,
\end{equation}
with $d=2$ here, complemented by initial conditions
\begin{equation*}
    \bm u(\bm x,0)=\bm u_0(\bm x), \quad \bm x\in \Omega,
\end{equation*}
and boundary conditions that will be described case-by-case. The solution $\bm u\in \R^m$ must stay in an open  subspace $\mathcal{D}\subset \R^d$ so that the flux $\bm f$ is defined. This is the convex invariant domain, which will be specified according to the problems studied.

For  scalar conservation laws with $\bm f$ linear taking the form of $\bm f(\bm u,\bm x)=\bm a(\bm x)\bm u$ or non-linear, $\bm u$ is a function with values in $\R$s but the solution must stay in $\big[\min\limits_{\bm x\in \R^d}\bm u_0(\bbx), \max\limits_{\bbx\in \R^d}\bm u_0(\bm x)\big]$, due to Kruzhkov's theory. In these particular examples, we set the invariant domain as $\mathcal{D}=\big[\min\limits_{\bm x\in \R^d}\bm u_0(\bm x), \max\limits_{\bm x\in \R^d}\bm u_0(\bm x)\big]$.

In the  case of the Euler equations of gas dynamics, $\bm u=(\rho, \rho \bm v, E)^T$, $\rho$ is the density, $\bm v=(u,v)^T$ is the velocity vector, and $E$ is the total energy, i.e., the sum of the internal energy $e$ and the kinetic energy $\tfrac{1}{2}\rho \bm v^2$. The flux is defined by
\begin{equation*}
    \bm f(\bm u)=\begin{pmatrix}
\rho \bm v\\
\rho \bm v\otimes \bm v+p\text{Id}_d\\
(E+p)\bm v
\end{pmatrix}
\end{equation*}
with $\text{Id}_d$ being a $d\times d$ identity matrix. We have introduced the pressure that is related to the internal energy and the density by an equation of state. Here we make the choice of a perfect gas,
\begin{equation*}
    p=(\gamma-1)e
\end{equation*}
with $\gamma$ denoting the specific heat ratio. In this example, the convex invariant domain is 
\begin{equation*}
    \mathcal{D}=\{\bm u=(\rho, \rho \bm v, E)^T\in \R^4, \text{such that}~\rho>0, e=E-\frac{1}{2}\rho\bm v^2>0\},
\end{equation*}
which is exactly equivalent to
\begin{equation}\label{Euler_GQL}
   \mathcal{D}_{\bm \nu}=\{\bm u=(\rho, \rho\bm v, E)^T\in\R^4 \text{ such that for all } \bm n_*\in \mathcal{N}, \bm u^T\bm n_*>0\}, 
\end{equation}
where, through the GQL (Geometric Quasilinearization) approach \cite{wu2023geometric}, the vector space $\mathcal{N}$ is given by
\begin{equation}\label{vector}
     \mathcal{N}=\{(1,0,0,0)^T\}\cup\{ (\frac{\Vert\bm\nu\Vert^2}{2} ,  -\bm\nu, 1)^T, \bm\nu\in \mathbb{R}^2\}.
\end{equation}
The reason is that the scalar product of $\bbu$ with 
$\big ( 1, \mathbf{0}, 0\big )^T$  is equal to $\rho$, and that of 
$\big ( -\frac{1}{2}\Vert \bbw\Vert , \bbw , 1\big )^T$ with $\bbu$ is equal to $e+\frac{1}{2}\Vert \bbv-\bbw\Vert^2.$

We note that for smooth solutions, \eqref{eq:conservative} can be re-written in a non-conservative form,
\begin{equation}\label{eq:non-conservative}
    \dpar{\bm u}{t}+\bm J\cdot\nabla \bm u=0,
\end{equation}
where
\begin{equation*}
    \bm J\cdot\nabla \bm u=A\dpar{\bm u}{x}+B\dpar{\bm u}{y},
\end{equation*}
with $A$ and $B$ being the Jacobians of the flux in the $x$- and $y$- direction, respectively. The system is hyperbolic, i.e., for any vector $\bm n=(n_x, n_y)$, the matrix
\begin{equation*}
    \bm J\cdot\bm n=A n_x+Bn_y
\end{equation*}
is diagonalisable in $\R$.

Here, we focus on the multidimensional extensions over unstructured triangular meshes in \cite{abgrall2024activefluxtriangularmeshes}.%,Liu_PAMPA_SW2D}.
 This new procedure, which is named as PAMPA (Point-Average-Moment PolynomiAl-interpreted) scheme, uses the degrees of freedom (DoFs) of the quadratic polynomials, which are all on the boundary of the elements, and the additional degree of freedom (DoF) that used is the cell average of the solution within each element. This leads to a potentially third-order method, with a globally continuous approximation. To evolve all DoFs in time, it relies on a standard Runge-Kutta time stepping method. To be specific, the cell average DoF is updated following the conservative formulation \eqref{eq:conservative}. Since the flux is assumed to be continuous crossing the cell interfaces, one can directly apply the divergence theorem and proper quadrature formulae to compute the fluctuations on each edge. The DoFs of the quadratic polynomials sitting on the edges are updated following the nonconservative formulation \eqref{eq:non-conservative}.  It was shown in \cite{Abgrall_camc, AB_HOAF} that the PAMPA scheme is linearly stable under a $1/2$ CFL condition. In order to prevent the numerical solution from producing spurious oscillations in the presence of discontinuities and possibly non-physical solutions. A stabilization technique equipped with the parachute low-order schemes and a MOOD mechanism \cite{Mood1} was applied in \cite{abgrall2024activefluxtriangularmeshes} to ameliorate the nonlinear instability issues or even code crashes. Another direction to ensure the nonlinear stability is, similar to those used in \cite{Kuzmin_MCL1,Zhang_MP,Guermond_IDP3, wu2023geometric}, to find a convex limiting approach between the high-order and low-order (parachute) fluxes, operators, or schemes so that the resulted method is invariant domain preserving (IDP). There are two attempts. One is in \cite{duan2024AF}, where the convex limiting method is applied to the average DoF and a scaling limiter with several flux vector splitting methods are used to the point value DoFs. One is in \cite{BPPAMPA1D}, where both the average DoF and point value DoFs are updated by a generic blending approach. To blend the fluxes and then ensure the updated average DoF stay in the convex invariant domain, it shares the similar philosophy as in most of the aforementioned IDP methods. For the update of point value DoFs on the boundaries of each element, the high-order residuals computed by a scheme based on Jacobian splitting are also blended with a low-order version computed by a first-order scheme with the local Lax-Friedrichs' flavor. This stem from the fact that no conservation requirement for evolving point value DoFs and a parachute IDP low-order scheme for point value DoFs does exist.  

In this paper, we extend the bound preserving method developed in \cite{BPPAMPA1D} to multidimensional case. The novelties of this work include the following: 
\begin{itemize}
    \item We propose a convex blending of the fluxes and residuals computed by the high-order PAMPA and low-order local Lax-Friedrichs schemes. With this blending, the updated solutions, including both average and point values, can be rewritten as a convex combination of previous time step solution and several intermediate solution states. This formulation facilitates the BP property analysis;
    \item We establish a sufficient condition for the BP property of the blended PAMPA scheme, which only requires the BP property of the intermediate solution states. This finding converts the goal of constructing BP PAMPA scheme to determine, through analysis, the optimal blending coefficients to reach the desired properties while trying to maintain as much as possible the high-order accuracy of the PAMPA scheme;
    \item We extend the BP PAMPA scheme to the
 Euler equations of gas dynamics. 
 The optimal blending coefficients can be directly obtained from the eigenvalues of a matrix with components only given by the intermediate solution states and the differences of high-order and low-order fluxes or residuals. This analysis is novel and accessible to multidimensional case; 
 \item  We implement the BP PAMPA scheme and demonstrate its robustness and effectiveness on several  numerical examples. %For instance, the KPP problem in the scalar conservation laws with non-convex flux, the LeBlanc, the double rarefaction Riemann, and the Sedov problems in the Euler equations of gas dynamics.  
\end{itemize}
 The time stepping will be done by a standard third order SSP Runge Kutta scheme, so we describe only the Euler forward step.
\section{Notations}
\subsection{Approximation space}
The approximation space is the same of the one adopted by the Virtual Finite Element Method (VEM) \cite{hitch}.  We first introduce some notations, following closely \cite{hitch}. The computational domain $\Omega$ is covered by {a set of} non-empty and non-overlapping polygons that are denoted by {$P$}. The notation $\bbx_P$ represents the centroid of {$P$}. {The elements} {$P$} may not be convex, but they are assumed to be star-shaped with respect to a point $\bby_P$ (that may be different of $\bbx_P$). For $D\subset \Omega$, the $L^2$ inner product between two functions in $L^2(D)$ is $\langle u,v\rangle_{D}$. When there is no ambiguity on $D$, we omit the subscript $D$. {For $\bsa=(\alpha_1, \alpha_2)$}, the scaled monomial of degree $\vert \bsa\vert=\alpha_1+\alpha_2$ are defined by
 \begin{equation}
 m_{\bsa}=\big ( \dfrac{\bbx-\bbx_P}{h_P}\big )^{\bsa}=\big ( \dfrac{x-x_P}{h_P}\big )^{\alpha_1} \big ( \dfrac{y-y_P}{h_P}\big )^{\alpha_2}, \qquad \bbx=(x,y).
 \end{equation}
 Here, $h_P$ is the diameter of $P$.%\footnote{Other choices are possible, such as a set of orthogonal polynomials \cite{orthVEM2017}. This set should contain the constant function $1$. This kind of choice allows for a better conditioning of the linear systems.}, according to the definition given in \cite{hitch}.
 The set of scaled monomial of degree $\vert\bsa\vert\leq k$ is a basis of $\P_k(P)$, which denotes the vector space of polynomials of degree less or equal to $k$, defined on ${P}$. In this paper, $k=2$. Similar definitions and comments can be done in three space dimensions. The scaled monomials are invariant by scaling: if $\hat{\bbx}=\lambda \bbx$, then 
\begin{equation}
\hat{m}(\hat{\bbx})=\big ( \dfrac{\hat{\bbx}-\bbx_{\hat{P} }}{h_{\hat{P}}}\big )^{\bsa}=m_{\bsa}=\big ( \dfrac{\bbx-\bbx_P}{h_P}\big )^{\bsa},
\end{equation}
because $\hat{P}$ is the image of $P$ by this mapping.
 
Now we introduce the degrees of freedom {of the approximation space} and define the local virtual space $V_2(P)$ for each $P$. A function $v_h\in V_2(P)$ is uniquely defined by the following setup:
 \begin{enumerate}
 \item $v_h$ is a polynomial of degree $2$ on each edge $e$ of $P$, that is $(v_h)_{|e}\in \P_2(e)$,
 \item $v_h$ is continuous on $\partial P$,
 \item $\Delta v_h$ is constant.
 \end{enumerate}
 Notice that $\P_2(P)\subset V_2(P)$, and a function of $V2k(P)$ is uniquely defined by the degrees of freedom given by:
 \begin{enumerate}
 \item The value of $v_h$ at the vertices of $P$,
 \item On each edge of $P$, the value of $v_h$ at the mid-points of  $3$  Gauss-Lobatto points on this edge (vertices and midpoint),
 \item The average value $m_{00}(v_h)$ of $v_h$ in $P$.
  \end{enumerate}
 The dimension of $V_2(P)$ is
 $  \dim V_2(P)=2{N_V} +1,$
  {with $N_V$ representing the number of vertices of $P$. The total number of degrees of freedom are then referred to with $N_{dofs}:=\dim V_2(P)$. Let $\{\varphi_i\}_{i=1}^{N_{dofs}}$ be the canonical basis for $V_2(P)$. We can then represent each $v_h \in V_k(P)$ in terms of its degrees of freedom by means of a Lagrange interpolation:
  	\begin{align}\label{eq:vh_repr} 
  		v_h = \sum_{i=1}^{N_{dofs}} \text{dof}_i(v_h)\varphi_i.
  	\end{align}
  For this basis, the usual interpolation property holds true:
  	\begin{equation}
  		\text{dof}_i(\varphi_j) = \delta_{ij}, \qquad i,j=1, \ldots, N_{dofs}.
  \end{equation}}
 The first step is to construct a projector $\pi^\nabla$ from $V_2({P})$ onto $\P_2({P})$. It is defined by two sets of properties. {First, for any $v_h\in V_2({P})$, the orthogonality condition
 \begin{equation}
 \label{projo:1}\big \langle \nabla p_2, \nabla \big ( \pi^\nabla v_h-v_h\big )\big \rangle =0, \quad \forall p_2\in \P_2({P}),
 \end{equation}
 has to hold true, which is defined up to the projection onto constants $P_0: V_k(P)\rightarrow \P_0(P)$, that can be fixed as follows:
   \begin{equation}
  P_0(v_h)=\frac{1}{\vert {P}\vert }\int_P v_h \; d\bbx =m_{(00)}(v_h).
  \end{equation}
 Then, we ask that 
\begin{equation}
 \label{projo:2}P_0( \pi^\nabla v_h-v_h)=0.
 \end{equation}
 We can \textit{explicitly} compute the projector by using only the degrees of freedom {previously introduced}, see \cite{hitch}, so that the exact shape of the basis functions is not  needed. This explains the adjective "virtual".

\section{Numerical schemes}
The computational domain $\Omega$ is covered by a family of nonoverlapping polygons $P$. The faces of the polygons will be generically denoted by $f$.  Because of lack of space, we skip the description of what we do on the boundaries, see \cite{BPPampa} for more details.

The numerical approximation at any time step is described by its average $\xbar{\bbu}_P^n$ for any polygon $P$, as well as the point value approximations $\bbu_\sigma^n$ for any Gauss-Lobatto points. There are two sets of updates, one for the average values and one for the point values.

%%%%%%%%%%%%%%%%%%%%%%%%%%%%%%
\section{Evolving point values}\label{sec:evoving_point_values}
%%%%%%%%%%%%%%%%%%%%%%%%%%%%%%
The point values will be evolved as follows:
\begin{subequations}\label{Pt:val}
\begin{equation}\label{Pt:val:1}
  u_\sigma^{n+1}=u_\sigma^n-\dt\sum_{P,\sigma\in P}{\Phi}_\sigma^P
\end{equation}
where
\begin{equation}\label{Pt:val:2}
  {\Phi}_\sigma^P={\Phi}_\sigma^{P,\,LO}+\theta_\sigma\Delta{\Phi}_\sigma^P
\end{equation}
\end{subequations}
 Following \cite{abgrall2024virtualfiniteelementhyperbolic},  we subdivide $P$ into sub-triangles as on Figure \ref{subtriangulation}--(a). The degrees of freedom on the boundary of $P$ are denoted by $\sigma_i$, where $i=0, \ldots , N_P-1$ with $N_P$ the number of DOFs on the boundary of $P$. The numbering is done clockwise modulo $N_P$, so that $\sigma_{-1}=\sigma_{N_P-1}$ and $\sigma_{N_P}=\sigma_0$. The vertices of the sub-triangles are $\{\sigma_i, \bby_P, \sigma_{i+1}\}$. The vertex $\sigma_i$ is shared by the triangles $\{\sigma_i, \bby_P, \sigma_{i+1}\}$ and $\{\sigma_{i-1},\bby_P,\sigma_{i}\}$.
The low-order residuals are defined as
\begin{equation}
  {\Phi}_{\sigma_i}^{E,\,LO}= \Phi_{\sigma_i}^{i-1,E,i}+
  \Phi_{\sigma_i}^{i,E,i+1}
\end{equation}
with
\begin{figure}[h]
\begin{center}
\subfigure[Sub-triangulation]{\includegraphics[width=0.35\textwidth]{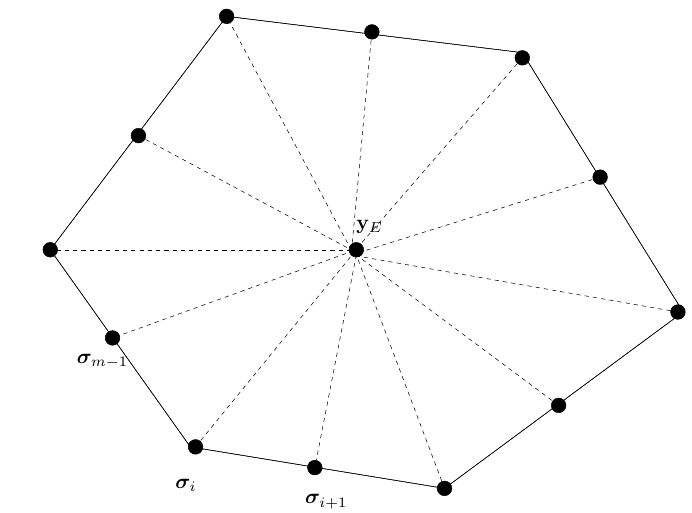}}\hspace*{0.2cm}
\subfigure[Normals]{\includegraphics[width=0.3\textwidth]{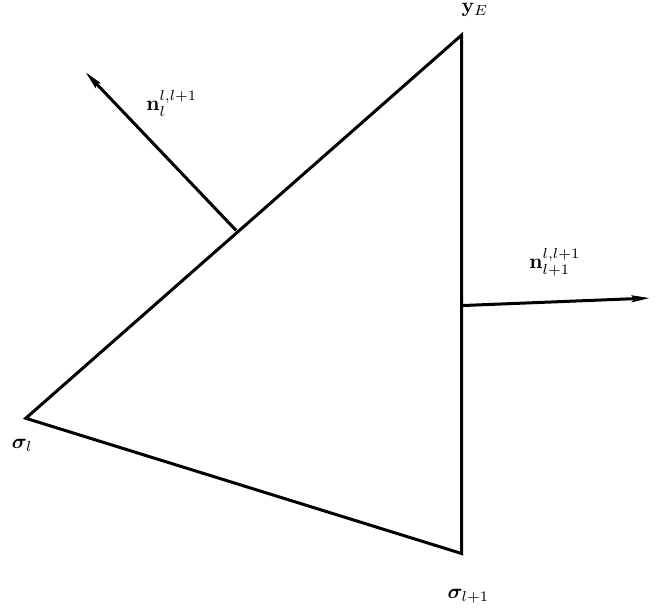}}
\end{center}
\caption{\label{subtriangulation} Sub-triangulation that is used to define a low order scheme. The point values on the boundary of $P$ are numbered clockwise modulo $N_P$ where $N_P$ is the number of DoFs on the boundary. The point $\bby_P$ is a point for which $P$ is star shaped, usually the centroid of $P$ in most applications. On the right panel are shown the normals.}
\end{figure}
\begin{equation}\label{res:LO}
\begin{split}
  \vert C_\sigma\vert\Phi_{\sigma_i}^{i-1,E,i}&=\frac{1}{6}\bigg ( \big (\bbf(\bbu_{\sigma_i})-\bbf(\xbar{\bbu}_E)\big )\cdot \bbn_i^{i-1,i}+
  \big ( \bbf(\bbu_{\sigma_{i-1}})-\bbf(\xbar{\bbu}_E)\big )\cdot \bbn_{i-1}^{i-1,i}\bigg )\\
  &\qquad +\frac{\alpha_{E}}{6}(\bbu_{\sigma}-\xbar{\bbu}_{i-1,i}^P)\\
  \vert C_\sigma\vert\Phi_{\sigma_i}^{i, i+1,E}&=\frac{1}{6}\bigg ( \big (\bbf(\bbu_{\sigma_i})-\bbf(\xbar{\bbu}_E)\big )\cdot \bbn_i^{i,i+1}+
  \big ( \bbf(\bbu_{\sigma_{i+1}})-\bbf(\xbar{\bbu}_E)\big )\cdot \bbn_{i+1}^{i, i+1 }\bigg )\\
  &\qquad +\frac{\alpha_{E}}{6}(\bbu_{\sigma}-\xbar{\bbu}_{i, i+1}^E)\\
  \xbar{\bbu}_{i-1,i}^E&=\frac{1}{6}\big ( \bbu_{\sigma_i}+\bbu_{\sigma_{i-1}}+\xbar{\bbu}_E), \qquad
  \xbar{\bbu}_{i,i+1}^E=\frac{1}{3}\big ( \bbu_{\sigma_i}+\bbu_{\sigma_{i+1}}+\xbar{\bbu}_E)
  \end{split}
\end{equation}
and the normals $\bbn_i^{i-1,i}$ is the outward normal of the edge $[\bby_P,\sigma_i]$ of the triangle $\{\sigma_{i-1}, \bby_P, \sigma_i\}$, etc, see figure \ref{subtriangulation}. The area $\vert C_\sigma \vert$ is 
$\vert C_{\sigma_i}\vert =\sum_{E, \sigma\in E} \tfrac{1}{3}\big (\vert [\sigma_i,\sigma_{i+1},\bby_E]\vert+\vert [\sigma_{i-1},\sigma_{i},\bby_E]\vert \big ).$
If we look more carefully at \eqref{res:LO}, the normals $\bbn_i^{i-1,i}$ and $\bbn_i^{i,i+1}$ are opposite to each other, so that 
\begin{equation}\label{res:LO:bon}
    \begin{split}
        \vert C_\sigma\vert \Phi_{\sigma_i}^E&=\frac{1}{6} \bigg ( \big (\bbf(\bbu_{\sigma_{i-1}})-\bbf(\xbar{\bbu}_E)\big )\cdot \bbn_{i-1}^{i-1,i}+\alpha_E\big ( \xbar{\bbu}_E-\bbu_{\sigma_i}\big )\bigg )\\
        &\quad + \frac{\alpha_E}{6}\big ( \xbar{\bbu}_E-\bbu_\sigma\big )\\
        & \qquad+\frac{1}{6}\bigg ( \big ( \big (\bbf(\bbu_{\sigma_{i+1}})-\bbf(\xbar{\bbu}_E)\big )\cdot \bbn_{i+1}^{i,i+1}+\alpha_E\big ( \xbar{\bbu}_E-\bbu_{\sigma_{i+1}}\big )\bigg )
    \end{split}
\end{equation}
Finally, looking at \eqref{res:LO}, we see that the scheme will be domain invariant preserving if
\begin{equation}\label{alpha_E}
\begin{split}
\alpha_E\geq \max\big (\alpha\big ( \bbu_{\sigma_{i-1}}, \bar\bbu_E, \bbn_{i-1}^{i-1,i}\big ), 
\alpha\big ( \bbu_{\sigma_{i+1}}, \bar\bbu_E, \bbn_{i+1}^{i, i+1}\big ) \bigg )
%\alpha_E&=\max\big (\rho(\overline{\nabla \bbf}_{\xbar{\bbu}_E, \bbu_{\sigma_i}}),
%                                    \rho(\overline{\nabla \bbf}_{\bbu_{\sigma_i},\xbar{\bbu}_E}), 
%                                    \rho(\overline{\nabla \bbf}_{\bbu_{\sigma_{i-1}},\xbar{\bbu}_E}), 
%                                    \rho(\overline{\nabla \bbf}_{\bbu_{\sigma_{i+1}},\xbar{\bbu}_E}))\\
%                                    &\qquad \times 
%\max\big( \Vert \bbn_i^{i-1,i}\Vert, \Vert \bbn_{i-1}^{i-1,i}\Vert, \Vert \bbn_{i}^{i, i+1, }\Vert, \Vert %\bbn_{i+1}^{i, i+1, }\Vert \big ) 
\end{split}
\end{equation}
where $\alpha(\bbu_1,\bbu_2,\bbn)$ is the maximum speed obtained from the Riemann problem between the states $\bbu_1$ and $\bbu_2$ in the direction $\bbn$ as evaluated in \cite{GuermondPopovFast}.
%where $\overline{\nabla \bbf}_{\bbu, \bbv}$ represents the Roe matrix between the states $\bbu$ and $\bbv$. The estimate \eqref{alpha_E} can easily be improved, see \cite{ShuPerthame}.\todo{check}

The high order residuals are evolved with
\begin{equation}\label{res:HO}
\Phi_\sigma^{E,H}=\bbN_\sigma \bbK^+_{\sigma,E}
\bbJ_{\sigma} \big (\nabla \pi^\nabla \bbu) +\alpha_P\sum_{\xi \text{dof  in }P} \text{dof}_\xi(\varphi_\xi) \text{dof}_{\xi}(\bbu-\pi^\bot \bbu)
\end{equation} 
\begin{figure}[h]
\begin{center}
\subfigure[]{\includegraphics[width=0.35\textwidth]{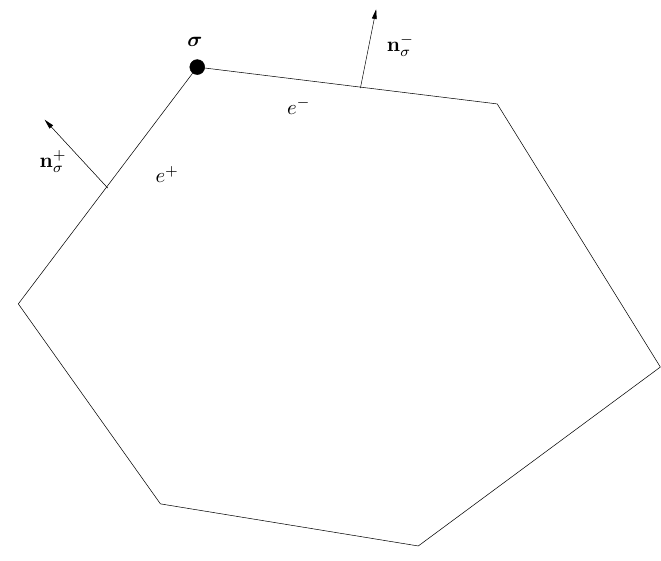}}
\subfigure[]{\includegraphics[width=0.35\textwidth]{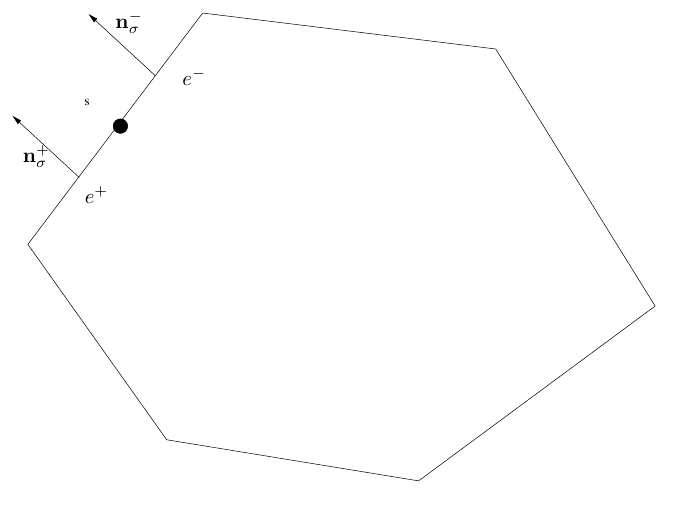}}
\end{center}
\caption{\label{normals} Definition of $e^\pm$ and $\bbn_\sigma^\pm$: (a) case of a vertex, (b) case of a non vertex. We use $\bbn_\sigma=\bbn_\sigma^++\bbn_\sigma^-$.}
\end{figure}where we have followed the following conventions (see \cite{abgrall2024virtualfiniteelementhyperbolic} and figure \ref{normals}):
\begin{equation*}
\begin{split}
\bbn_\sigma=\bbn_\sigma^++\bbn_\sigma^-,  & \qquad \bbJ_{\sigma}=\bbJ(\bbu_\sigma)\\
\bbK^+_{\sigma,E}=\bbJ_\sigma\cdot \bbn_{\sigma}, & \qquad \bbN_\sigma^{-1}=\sum_{E, \sigma\in E}\bbK^+_{\sigma,E}
\end{split}
\end{equation*}
and $\alpha_P$ is an upper bound of the eigen-speeds of $\mathbf{J}$. 
It is shown in \cite{abgrall2024virtualfiniteelementhyperbolic} that $\bbN_\sigma K_\sigma^+$ is always defined. The sum in \eqref{res:HO} is a form of dissipation inspired by the VEM litterature, it does not spoil the formal accuracy. The scalar $\alpha_P$ is an estimation of the maximum wave speed in $P$.

%Since the domain is bounded, we will also need boundary conditions. These are described in Section \ref{sec:BCs}. In order to maintain generality, we will denote by $\Phi_\sigma^{g,H}$ and $\Phi_\sigma^{f,L}$ the residuals constructed for $\sigma$ and the face $g\subset \Gamma$ to describe the boundary conditions.

 Following \cite{Guermond_IDP3,Kuzmin_MCL1,duan2024AF,BPPAMPA1D},  we set $\Phi_\sigma^E=\theta_\sigma^E\Phi_\sigma^{H,E}+(1-\theta_\sigma^E)\Phi_{\sigma}^{L,E}$ 
  so that 
the Euler forward update writes, where $\vert E_\sigma\vert$  and $\vert g\vert_\sigma$ are  left free for now, as:
\begin{equation}\label{scheme:residual}
\begin{split}
\bbu_\sigma^{n+1}&=\bbu_\sigma^n-\dt \sum_{E, \sigma \in E} \Phi_\sigma^E(\bbu)\\
%\Delta t \sum_{g\subset \Gamma, \sigma\in g}\Phi_{\sigma}^d(\bbu)\\
&=\bbu_\sigma^n-\dt \sum_{E, \sigma \in  E} \bigg ( \Phi_\sigma^E(\bbu)+\frac{\alpha_E}{|E|_\sigma}\big ( \bbu_\sigma^n-\bbu_\sigma^n)\bigg )
%-\Delta t\sum\limits_{g\subset \Gamma, \sigma\in g} 
%\bigg ( \Phi_\sigma^g(\bbu)+\frac{\alpha_g}{\vert g\vert_\sigma}\big (\bbu_\sigma^n-\bbu_\sigma^n\big )\bigg )
\\
&=\bigg (1-\dt \big (\sum_{E, \sigma\in E} \frac{\alpha_E}{|E|_\sigma}\big )\bigg ) \bbu_\sigma^n+
\sum_{E, \sigma\in E}\dt\frac{\alpha_E}{\vert E_\sigma\vert} \bigg ( \bbu_\sigma^n-\vert E_\sigma\vert\frac{\Phi_\sigma^E}{\alpha_E}\bigg )
%+
%\sum_{g\subset\Gamma, \sigma\in g}\Delta t \frac{\alpha_g}{\vert g\vert_\sigma}\bigg ( \bbu_\sigma^n-\vert g\vert_\sigma\frac{\Phi_\sigma^g}{\alpha_g}\bigg )
\\
&=\bigg (1-\dt \big (\sum_{E, \sigma\in E} \frac{\alpha_E}{|E|_\sigma}\big )\bigg ) \bbu_\sigma^n+\sum_{E, \sigma\in E}\dt\frac{\alpha_E}{\vert E_\sigma\vert} \tilde{u}_\sigma^E
\end{split}
\end{equation}
with
\begin{equation}\label{scheme:residual:star}
\begin{split}
\tilde{u}_\sigma^E&= \bbu_\sigma^n-\vert E_\sigma\vert\frac{\Phi_\sigma^E}{\alpha_E}=
\bbu_\sigma^n-\vert E_\sigma\vert \frac{\Phi_\sigma^{L,E}}{\alpha_E}-\theta_\sigma^E\vert E_\sigma\vert \frac{\Delta \Phi_\sigma^E}{\alpha_E}\\
&=\bbu_\sigma^{E,\star}-\theta_\sigma^E \vert E_\sigma\vert\frac{\Delta \Phi_\sigma^E}{\alpha_E}.
\end{split}
\end{equation}
%and similarly for $\widetilde{\bbu}_\sigma^g$,
%\begin{equation}\label{scheme:residual:star:bc}
%   \widetilde{\bbu}_\sigma^g=\bbu_{\sigma}^{g,\star}-\theta_\sigma^g\vert g\vert_\sigma\frac{\Delta \Phi_\sigma^g}{\alpha_g} \text{ and }  \bbu_{\sigma}^{g,\star}=\bbu_\sigma^n-\vert g\vert \frac{\Delta \Phi_{\sigma}^{L, g}}{\alpha_g}
%\end{equation}
If $\bbu_\sigma^n\in \mathcal{D}$ and $\bbu_\sigma^{g,\star}\in \DD$, 
\begin{equation}\label{cfl:1}
\dt \sum_{E, \sigma\in E} \frac{\alpha_E}{|E|_\sigma}\leq 1,
\end{equation}
if $\bbu_\sigma^{E,\star}\in \mathcal{D}$ and 
\begin{equation}\label{theta:pt}
\begin{split}
\theta_\sigma^E&=\min\bigg (1,\min\limits_{E, \sigma\in E} \frac{\alpha_E}{\vert E_\sigma\vert}\dfrac{\min\big (M- \bbu_\sigma^{E,\star},
\bbu_\sigma^{E,\star}-m\big )}{\Delta \Phi_\sigma^E}\bigg ), 
\end{split}
\end{equation}
then $\bbu_\sigma^{n+1}\in \mathcal{D}$.

In order to get $\bbu_\sigma^{E,\star}\in \mathcal{D}$, we need to find out what are $\alpha_E$, $\alpha_g$, $\vert E_\sigma\vert $, and $\vert g_\sigma\vert$. We do this for $\alpha_E$ and $\vert E_\sigma\vert $. 
Using \eqref{res:LO:bon}, we write $\bbu_\sigma^{E,\star}$ as 
\begin{equation*}
    \begin{split}
        \bbu_\sigma^{E,\star}&=\bbu_\sigma^n-\frac{\vert E_\sigma\vert }{\alpha_E}\Phi_\sigma^E\\
        &=\frac{1}{3}\bigg ( \xbar{\bbu}_E^n-\frac{\vert E_\sigma}{2\alpha_E\vert C_\sigma\vert} \bigg ( \big (\bbf(\bbu_{\sigma_{i-1}})-\bbf(\xbar{\bbu}_E)\big )\cdot \bbn_{i-1}^{i-1,i}+\alpha_E\big ( \xbar{\bbu}_E-\bbu_{\sigma_i}\big )\bigg )\bigg )\\
        &\quad +\frac{1}{3}\bigg ( \xbar{\bbu}_E^n- \frac{\vert E_\sigma}{2\alpha_E\vert C_\sigma\vert}\big ( \xbar{\bbu}_E-\bbu_\sigma\big )\bigg )\\
        &\qquad +\frac{1}{3}\bigg (\xbar{\bbu}_E^n-\frac{\vert E_\sigma}{2\alpha_E\vert C_\sigma\vert}\bigg ( \big ( \big (\bbf(\bbu_{\sigma_{i+1}})-\bbf(\xbar{\bbu}_E)\big )\cdot \bbn_{i+1}^{i,i+1}+\alpha_E\big ( \xbar{\bbu}_E-\bbu_{\sigma_{i+1}}\big )\bigg )
    \end{split}
\end{equation*}
This suggest to set 
$\vert E_\sigma\vert=\vert C_\sigma\vert$%,\label{choiceE}
and then we get
\begin{equation*}
    \begin{split}
        \bbu_\sigma^{E,\star}&=\bbu_\sigma^n-\frac{\vert E_\sigma\vert }{\alpha_E}\Phi_\sigma^E\\
        &=\frac{1}{3}\bigg ( \frac{\xbar{\bbu}_E^n+\bbu_{\sigma_{i-1}}^n}{2}-\frac{1}{2\alpha_E}  \big (\bbf(\bbu_{\sigma_{i-1}})-\bbf(\xbar{\bbu}_E)\big )\cdot \bbn_{i-1}^{i-1,i}\bigg )\\
        &\quad +\frac{1}{3} \xbar{\bbu}_E^n 
        +\frac{1}{3}\bigg (\frac{\xbar{\bbu}_E^n+\bbu_{\sigma_{i+1}}^n}{2}-
        \frac{ \big (\bbf(\bbu_{\sigma_{i+1}})-\bbf(\xbar{\bbu}_E)\big )\cdot \bbn_{i+1}^{i,i+1}}{2\alpha_E} \bigg )
    \end{split}
\end{equation*}
This shows that with this choice, $\bbu_\sigma^\star\in \DD$ under the condition
\begin{equation}\label{cfl:pt}
\Delta t \frac{\sum\limits_{E, \sigma\in E}\alpha_E}{\vert C_\sigma\vert}\leq 1.
\end{equation}

 %%%%%%%%%%%%%%%%%%%%%%%%%%%%%%%%
\section{Evolving  cell averages}\label{sec:evolve_cell_average}
%%%%%%%%%%%%%%%%%%%%%%%%%%%%%%%%%
 In the following,  $\FF_E$ is the set of faces interior to $\Omega$.
The cell averages will be updated as follows:
\begin{subequations}\label{ave:scheme}
\begin{equation}\label{ave:scheme:1}
    \xbar{\bbu}_E^{n+1}=\xbar{\bbu}_E^n-\frac{\dt}{\vert E\vert} \sum_{f\in \FF_E }\vert f\vert \; \hbbf(\bbu_E,\bbu_{E_f}, \bbn_f).
\end{equation}
Here$\vert f\vert$ is the measure of  $f$, $E_f$ is the polygon on the other side of $f$, $\bbn_f$ the outward normal vector. The flux $\hbbf$ is written as
\begin{equation}\label{ave:scheme:2}
\hbbf(\bbu_E,\bbu_{E_f},\bbn_f)=\hbbf(\bbu_E,\bbu_{E_f},\bbn_f)^{LO}+\eta_{f}\Delta \hbbf(\bbu_E,\bbu_{E_f},\bbn_f).
\end{equation}
\end{subequations}

The low-order numerical flux is defined as
\begin{equation}\label{flu:LO}
\hbbf(\bbu_E,\bbu_{E_f},\bbn_f)^{LO}=\frac{\bbf(\xbar{\bbu}_E)\cdot \bbn_f+\bbf(\xbar{\bbu}_{E_f})\cdot\bbn_f}{2}-\frac{\beta_l}{2}(\xbar{\bbu}_{{E_f}}-\xbar{\bbu}_E),
\end{equation}
while the high order flux is
\begin{equation*}
\hbbf_f^{HO}(\bbu_E,\bbn)=\int_{f} \bbf(\bbu)\cdot \bbn_f \; d\gamma.
\end{equation*} 
This is evaluated by quadrature formula, and here we exploit the fact that the point values are localized at Gauss-Lobatto points,
\begin{equation}\label{flu:HO}
\hbbf_f^{HO}(\bbu_E,\bbn)=\vert f\vert \sum_{\sigma\in f} \omega_{\sigma} \bbf(\bbu_\sigma)\cdot \bbn_f.
\end{equation} 

Similarly to the point values evaluation, the Euler forward step can be rewritten as 
%\begin{equation*}
%\begin{split}
%\xbar{\bbu}_E^{n+1}&=\xbar{\bbu}_E^n-\frac{\dt}{\vert E\vert }\sum_{f\in \FF_E}\vert f\vert 
% \hbbf(\xbar{\bbu}_E^n, \xbar{\bbu}_{E_f}^n, \bbn_f)\\
% &=\xbar{\bbu}_E^n-\frac{\dt}{\vert E\vert }\sum_{f\in \FF_E}\vert f\vert 
% \hbbf(\xbar{\bbu}_E^n, \xbar{\bbu}_{E_f}^n, \bbn_f)-\frac{\dt}{\vert E\vert} \sum_{f\in \FF_E} \vert f\vert\big ( \alpha_{E_f}\xbar{\bbu}_E^n - \alpha_{E_f}\xbar{\bbu}_{E}^n+\bbf(\xbar{\bbu}_E)^n\cdot \bbn_f \big )\\
%% &\qquad+
%% \frac{\dt}{\vert E\vert }\sum_{f\in \FF^b_E}\vert f\vert 
%% \hbbf(\xbar{\bbu}_E^n \bbn_f)-\frac{\dt}{\vert E\vert} \sum_{f\in \FF^b_E} \vert f\vert\big ( \alpha_{E_b}\xbar{\bbu}_E^n - \alpha_{E_b}\xbar{\bbu}_{E}^n+\bbf(\xbar{\bbu}_E)^n\cdot \bbn_f \big )\\
% %
% &=\bigg ( 1-\frac{\dt}{\vert E\vert }\sum_{f\in \FF^E}\vert f\vert \alpha_{{E_f}}\bigg ) \xbar{\bbu}_E^n+
% \sum_{f\in \FF_E} \dt \frac{\vert f\vert \alpha_{{E_f}}}{\vert E\vert } \bigg ( \xbar{\bbu}_E^n-\dfrac{ \bbf(\xbar{\bbu}_E^n)\cdot \bbn_f-\hbbf (\xbar{\bbu}_E^n,\xbar{\bbu}_{{E_f}}^n,\bbn_f)}{\alpha_{{E_l}}}\bigg )\\
%% &\qquad \qquad + 
%% \sum_{f\in \FF^b_E} \dt \frac{\vert f\vert \alpha_{{E_b}}}{\vert E\vert } \bigg ( \xbar{\bbu}_E^n-\dfrac{ \bbf(\xbar{\bbu}_E^n)\cdot \bbn_f-\hbbf^b(\xbar{\bbu}_E^n,\bbn_f)}{\alpha_{{E_b}}}\bigg )
%\end{split}
%\end{equation*}
%then, using \eqref{ave:scheme:2},% and \eqref{ave:scheme:b},
%we get
\begin{equation*}
    \begin{split}
        \xbar{\bbu}_E^{n+1}=\bigg ( 1-\frac{\dt}{\vert E\vert }\sum_{f\in \FF^E}\vert f\vert \alpha_{{E_f}}
        \bigg ) \xbar{\bbu}_E^n+\sum_{f\in \FF_E}\Delta t \frac{\vert f\vert}{\vert E\vert} \widetilde{\xbar{\bbu}}_E    \end{split}
\end{equation*}
where
$$
\widetilde{\xbar{\bbu}}_E=\xbar{\bbu}_E^n-\dfrac{ \bbf(\xbar{\bbu}_E^n)\cdot \bbn_l-\hbbf (\xbar{\bbu}_E^n,\xbar{\bbu}_{{E_l}}^n,\bbn_l)}{\alpha_{{E_l}}}
=\xbar{\bbu}^{f_l, E \star}-\theta_{f_l}\dfrac{\Delta \hbbf(\xbar{\bbu}_E^n,\xbar{\bbu}_{{E_l}}^n,\bbn_l)}{\alpha_{{E_l}}}
$$
with
$$\xbar{\bbu}^{f_l,  \star}=\frac{\xbar{\bbu}_E^n+\xbar{\bbu}_{{E_l}}^n}{2}-\frac{\bbf(\xbar{\bbu}_E^n)\cdot \bbn_l-\bbf(\xbar{\bbu}_{{E_l}}^n)\cdot \bbn_l}{2\alpha_E}.$$

In the end, we get that if $\xbar{\bbu}_P\in \DD$ for $P=E$ and all the other polygons sharing a face with $E$,  if 
\begin{subequations}\label{flux:positif}
    \begin{equation}
        \label{flux:positif:alpha}
\alpha_E\geq \alpha(\xbar{\bbu}_E^n),\xbar{\bbu}_{{E_f}}^n), \bbn_f \big ),
\end{equation}
\begin{equation}\label{flux:positif:cfl}
\frac{\dt}{\vert E\vert } \sum_{f\in \FF_E}\vert f\vert \alpha_{{E_f}}\leq 1
\end{equation}and
\begin{equation}\label{flux:positif:theta}\theta_{f}=\min\big (1,\dfrac{\min(M-\xbar{\bbu}^{f_l,  \star}, \xbar{\bbu}^{f_l,  \star}-m)}{\vert \Delta \hbbf(\bbu_E, \bbu_{{E_l}}, \bbn)\vert}\big ), \end{equation}
\end{subequations}
then $\xbar{\bbu}_E^{n+1}\in \DD.$ 

The final condition on the time step is
\begin{equation}\label{cfl:ave}
    \Delta t \bigg ( \frac{\sum_{f\in \FF_E} \vert f\vert \alpha_{E_l}
     }{\vert E\vert }\bigg )\leq 1.
\end{equation}
\section{System: the example of the Euler equations}
%Detailler la valeur optimale de theta
The Euler equations in $d$ dimensions are
\begin{equation}
\label{eq:euler}
\dpar{\bbu}{t}+\text{ div }\bbf(\bbu)=0
\end{equation}
where
$$\bbu=\begin{pmatrix} \rho\\ \rho \bbv \\ E\end{pmatrix}, \quad
\bbf(\bbu)=\begin{pmatrix}
\rho\bbv \\
\rho \bbv\otimes \bbv+p\Id\\
(E+p)\bbv
\end{pmatrix}$$
The density is $\rho$, the velocity is $\bbv$, $E=e+\frac{1}{2}\rho \bbv^2$ is the total energy, $e$ the internal one and the pressure is $p=p(\rho, e)$. We will consider here the simplest case of a calorically perfect gas, where
$$p=(\gamma-1)e=(\gamma-1)\big (E-\frac{1}{2}\rho \bbv^2\big ).$$
The system is symetrizable, with the entropy defined by $S=\rho \big ( \log p-\gamma\log \rho\big )$, and hence the system is hyperbolic.

The system \eqref{eq:euler} is well defined as long as the density and the internal energy are positive, i..e the invariant domain is
$$\DD=\{\bbu, \text{ such that } \rho\geq 0, e\geq 0\}.$$
Here we will use the geometric characterisation defined by \eqref{Euler_GQL}.

The scheme is defined by \eqref{Pt:val} and \eqref{ave:scheme} with flux are residuals defined by \eqref{Pt:val:2} and \eqref{ave:scheme:2}, the high-order flux and residuals defined by \eqref{res:LO}, \eqref{res:HO}, \eqref{flu:LO}, \eqref{flu:HO}. Of course the conserved variable and the flux are replaced by the Eulerian ones. The only thing to define are the blending parameters.  To simplify the text, we again avoid to consider boundary conditions. 

To guaranty the positivity of the density, it is enough to take
$$\theta_{f_l}\geq \min(1,  \alpha_{f_l}\dfrac{
\rho^{E,\star}}{\vert \Delta \hbbf_\rho^{f_l}\vert}):=\widetilde{\theta^\rho_{f_l}}, \quad\theta_\sigma\geq \min(1,  \alpha_E\dfrac{ \rho_\sigma^{E,\star}}{\vert \Delta \Phi_{\sigma,\rho}^E\vert})
:=\widetilde{\theta^\rho_{\sigma}}.$$

To guaranty the positivity of the internal energy, and setting for $\bbw\in \R^d$, 
$\psi(\bbw)=(\frac{\Vert \bbw\Vert^2}{2}, -\bbw, 1)$,  we impose in addition
$\theta_{f_l}\geq 
\min(1,  \alpha_{f_l}\min\limits_{\bbw\in \R^d}\dfrac{\bbu^{f_l,\star}\psi(\bbw)^T}{\vert \big (\Delta \hbbf^{f_l}\big )^T\psi(\bbw)\vert}):=\widetilde{\theta^e_{f_l}}$
 \text{ and } $\theta_\sigma\geq 
 \min(1,  \alpha_{E'}\min\limits_{\bbw\in \R^d}\dfrac{\bbu^{E,\star}\psi(\bbw)^T
}{\vert \big (\Delta \Phi_\sigma^E\big )^T\psi(\bbw)\vert})
:=\widetilde{\theta^e_{\sigma}},$
so that we take
$$\theta_{f_l}=\min\big (\widetilde{\theta^\rho_{f_l}},\widetilde{\theta^e_{f_l}}\big )  \text{ and } \theta_\sigma=\min \big (\widetilde{\theta^\rho_{\sigma}},\widetilde{\theta^e_{\sigma}}\big ).$$
Of course the issue is to evaluate
$$\min\limits_{\bbw\in \R^d}\dfrac{\bbu^{E,\star}\psi(\bbw)^T
}{\vert \big (\Delta \Phi_\sigma^E\big )^T\psi(\bbw)\vert} \text{ and } 
\min\limits_{\bbw\in \R^d}\dfrac{\bbu^{f_l,\star}\psi(\bbw)^T
}{\vert \big (\Delta \hbbf^{f_l}\big )^T\psi(\bbw)\vert}. $$
We discuss this for the flux, it is exactly the same for the residual.  The first thing to notice is that
$$\min\limits_{\bbw\in \R^d}\dfrac{\bbu^{f_l,\star}\psi(\bbw)^T
}{\vert \big (\Delta \hbbf^{f_l}\big )^T\psi(\bbw)\vert}=\bigg ( \max\limits_{\bbw\in \R^d}\dfrac
{\vert \big (\Delta \hbbf^{f_l}\big )^T\psi(\bbw)\vert}
{\bbu^{f_l,\star}\psi(\bbw)^T}
\bigg )^{-1}, $$
then that
$$\dfrac
{\vert \big (\Delta \hbbf^{f_l}\big )^T\psi(\bbw)\vert}
{\bbu^{f_l,\star}\psi(\bbw)^T}=\dfrac{\vert  \frac{\alpha_1}{2}\Vert \bbw\Vert^2-\bba_1^T\bbw+\beta_1\vert }
{ \frac{\alpha_0}{2}\Vert \bbw\Vert^2-\bba_0^T\bbw+\beta_0}
$$
Replacing $\bbw$ by $\bbw/w_{d+1}$ because we can choose any $\bbw$, we get
$$\dfrac
{\vert \big (\Delta \hbbf^{f_l}\big )^T\psi(\bbw/w_{d+1})\vert}
{\bbu^{f_l,\star}\psi(\bbw/w_{d+1})^T}=\dfrac{\vert  \frac{\alpha_1}{2}\Vert \bbw\Vert^2-\bba_1^T\bbw w_{d+1}+\beta_1w_{d+1}^2\vert }
{ \frac{\alpha_0}{2}\Vert \bbw\Vert^2-\bba_0^T\bbw w_{d+1}+\beta_0w_{d+1}^2}$$
because $w_{d+1}^2\geq 0$, and by assumption we will have
$$\frac{\alpha_0}{2}\Vert \bbw\Vert^2-\bba_0^T\bbw w_{d+1}+\beta_0w_{d+1}^2\geq 0$$ because again $w_{d+1}^2\geq 0$. Hence, setting $\bbz=(\bbw,w_{d+1})\in \R^{d+1}$ the problem reduces to study the Rayleigh quotient
$$\max\limits_{\bbz\in \R^{d+1}}\dfrac{\bbz^T B\bbz}{\vert \bbz^T B\bbz\vert }$$ which is simply
$\rho(A^{-1/2}BA^{-1/2})$ since, by assumption, saying that 
$\frac{\alpha_0}{2}\Vert \bbw\Vert^2-\bba_0^T\bbw w_{d+1}+\beta_0w_{d+1}^2>0$ for non zero $\bbz$ amounts to say that the symmetric matrix $A$ is positive definite.

The evaluation of $\rho(A^{-1/2}BA^{-1/2})$ can be done by some iterative method. It turns out that in the particular case we consider, we can find a simple analytical formula.
We have
$$A=\begin{pmatrix}
\alpha_0\Id_d & -2 \bba\\
-2\bba^T & 2\alpha_{d+1}\end{pmatrix}, \quad B=\begin{pmatrix}
\beta_0\Id_d & -2 \bbb\\
-2\bbb^T & 2\beta_{d+1}\end{pmatrix}$$ and after some easy calculations, we find that 
$\rho\big (A^{-1/2}BA^{-1/2}\big )=\max\big ( \frac{\vert\beta_0\vert}{\alpha_0}, \vert \lambda_+\vert, \vert\lambda_-\vert \big )$ with
$$\lambda_{\pm}=\dfrac{\big (
-4 \bba^T\bbb+\big (\beta_{d+1}\alpha_0+\beta_0\alpha_{d+1}\big )\big )\pm \sqrt{\Delta}
}{4\Vert \bba\Vert^2-2\alpha_0\alpha_{d+1}}$$
\text{ and  }$$\Delta=\big (4 \bba^T\bbb-\big (\beta_{d+1}\alpha_0+\beta_0\alpha_{d+1}\big )\big )^2-4\big (2\Vert \bba\Vert^2-\alpha_0\alpha_{d+1}\big )\big (2 \Vert \bbb\Vert^2-\beta_0\beta_{d+1}\big )\geq 0.$$

\begin{remark} 
This technique is not specific to the schemes developed in this paper but has a larger potential. It has also been used in \cite{GauthierLBM}, for a completely different scheme that uses a kinetic formulation of the Euler equations.
\end{remark}

\section{Numerical results}
All the solutions we presents have discontinuous features. The algorithm is exactly the same as that of \cite{abgrall2024virtualfiniteelementhyperbolic} for smooth solution and we refer to this paper for them (vortex and acoustic problem).
\subsection{Scalar problems: Zalesak test case}
The initial condition of this very classical test case of linear convection with circular rotation, $u_0$ defined in $-1,1]^2$ is $u_0(\bbx)=0$ except
\begin{equation*}
\begin{split}
\text{if }r^2=\Vert \bbx-(0.25,0.5)\Vert^2\leq 0.15& u_0(\bbx)=0.25(1+\cos(\pi\frac{r^2}{0.15})\\
\text{if }r^2=\Vert\bbx-(0.5,0.25)\Vert^2\leq 0.15& u_0(\bbx)=1-\frac{r^2}{0.15}\\
\text{if } r^2=\Vert\bbx-(0.5,0.75)\Vert^2\leq 0.15 & u_0(\bbx)=1.
\end{split}
\end{equation*}
The solution is presented after one rotation in figure \ref{fig:Zalesak}. Here, we have taken $m=0$ and $M=1$, the solution stays within these bounds. Compared to the exact solution, the results are excellent.
\begin{figure}[h!]
\begin{center}
\subfigure[]{\includegraphics[width=0.45\textwidth,clip=]{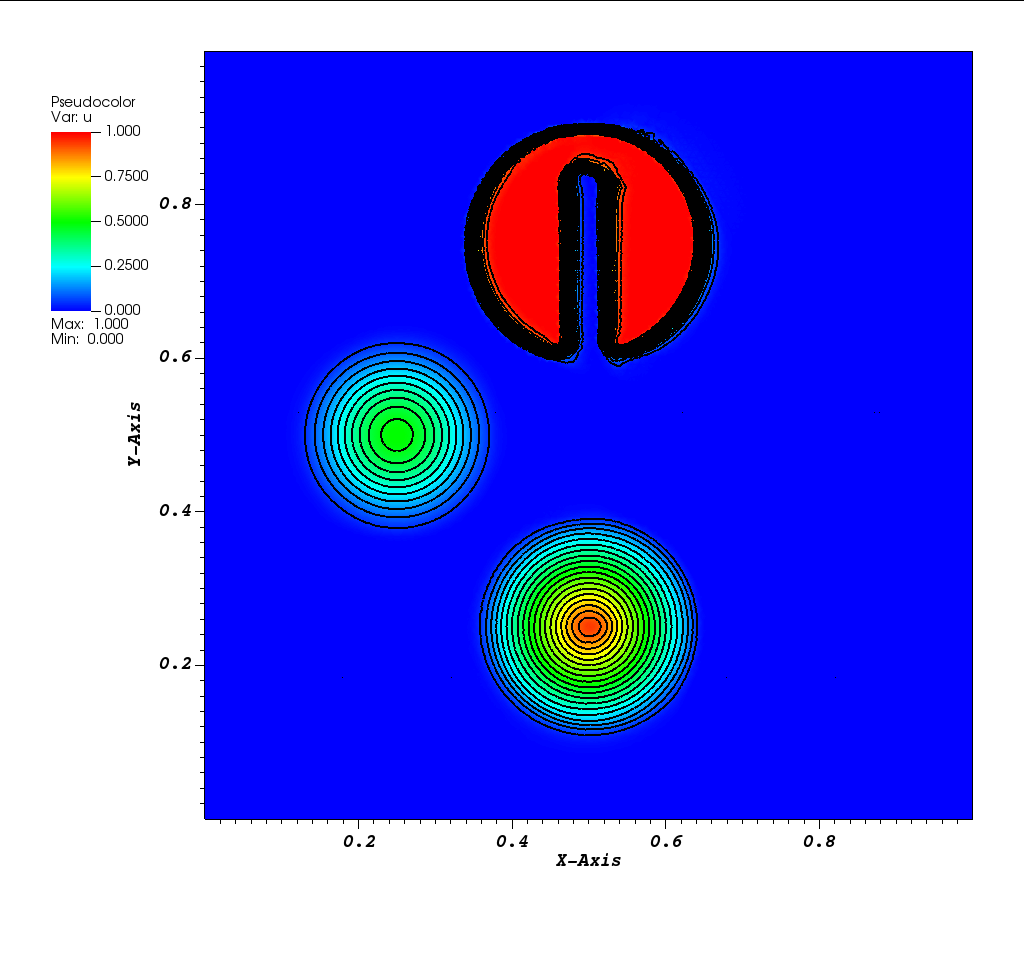}}
\subfigure[]{\includegraphics[width=0.45\textwidth,clip=]{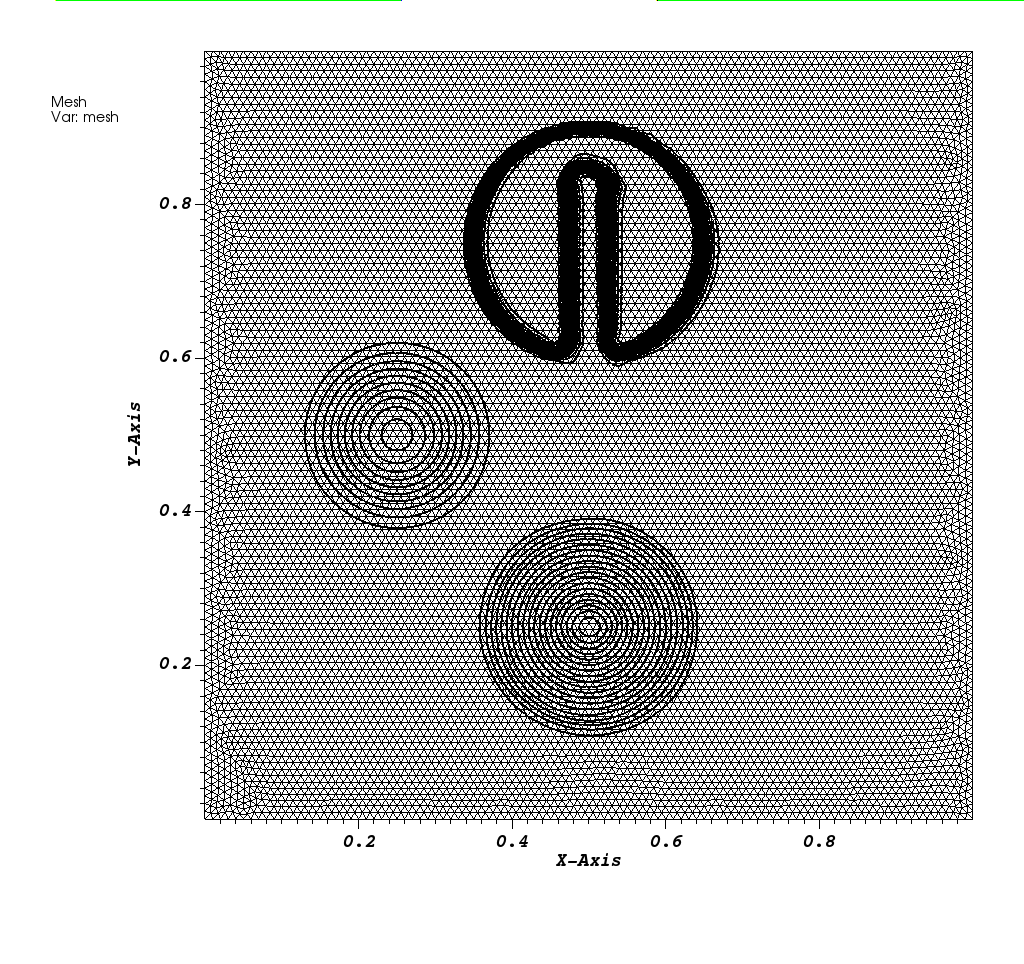}}
\subfigure[]{\includegraphics[width=0.45\textwidth,clip=]{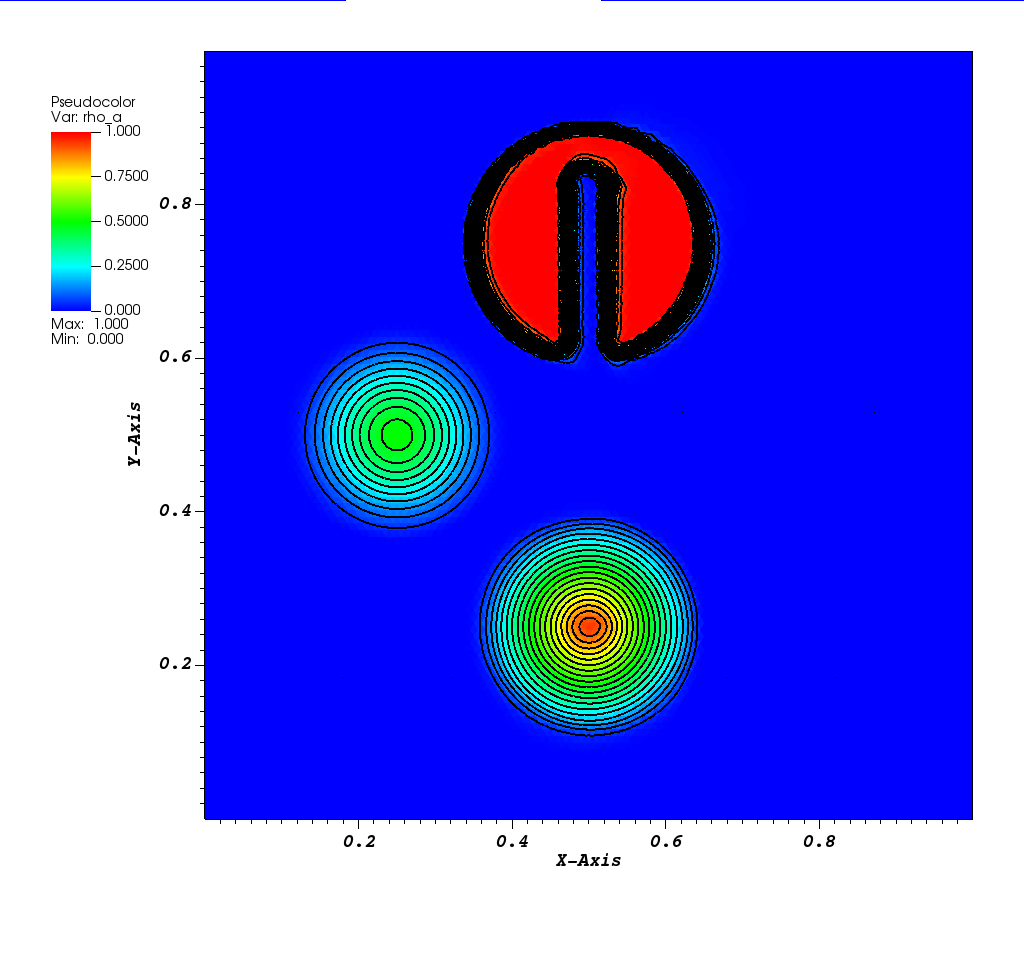}}
\subfigure[]{\includegraphics[width=0.45\textwidth,clip=]{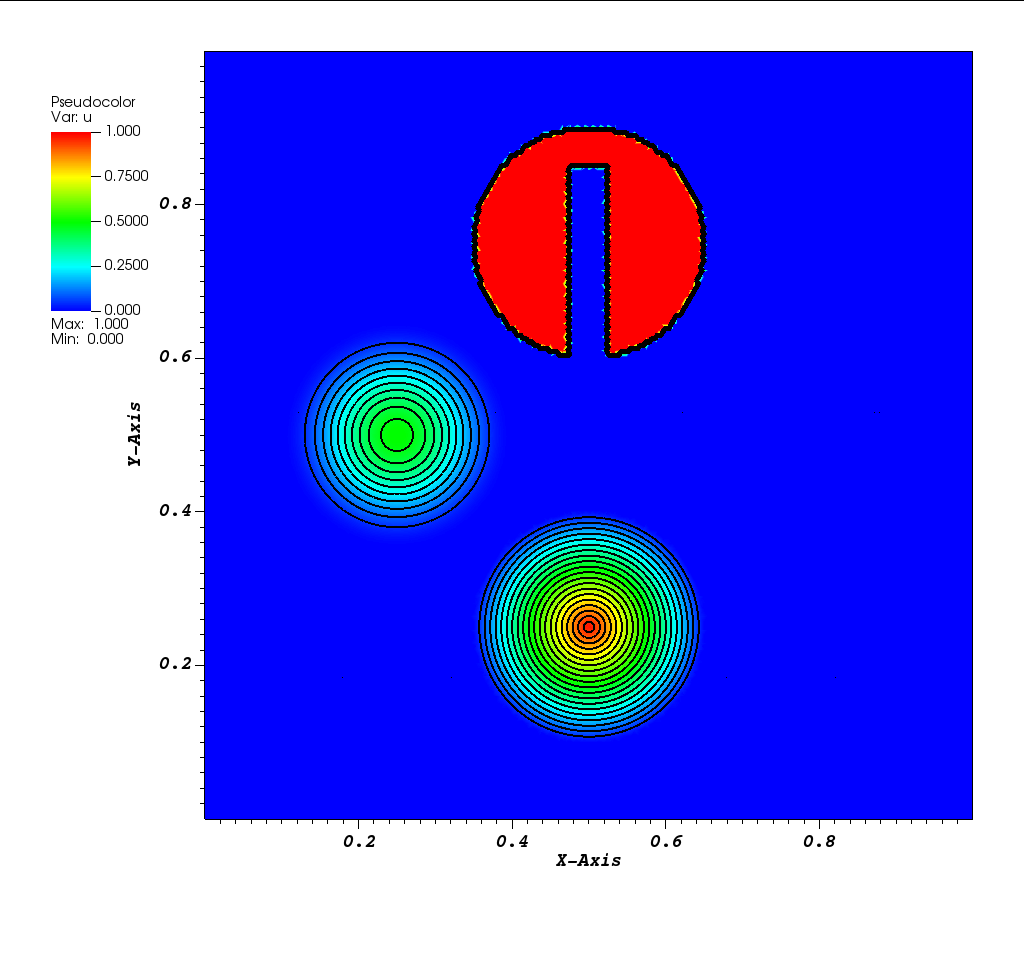}}
\caption{\label{fig:Zalesak} Zalesak test case, after 1 rotation. 20 isolines of (a): the point values  (b):the average, and (d): the exact solution.  CFL=0.4, $4709$ nodes and $23704$ elements. On (c) is represented the mesh and the average isolines}
\end{center}
\end{figure}
\subsection{KT case}
The initial condition is 
$$ (\rho, u,v,p)=\left \{\begin{array}{ll}
(\rho_1,u_1,v_1,p_1)=(1.5, 0, 0, 1.5)& \text{ if } x\geq1\text{ and } y\geq 1,\\
(\rho_2,u_2,v_2,p_2)=(0.5323, 1.206, 0, 0.3) & \text{ if } x\leq 1 \text{ and } y\geq 1,\\
(\rho_3,u_3,v_3,p_3)=(0.138, 1.206, 1.206, 0.029)&\text{ if } x\leq 1\text{ and }y\leq 1,\\
(\rho_4,u_4,v_4,p_4)=(0.5323, 0, 1.206, 0.3) &\text{ if } x\leq1\text{ and } y\leq 1.
\end{array}\right .
$$
Here, the four states are separated by shocks. The domain is $[-2,2]^2$. The solution at {$t_f=3$} is displayed in Figure \ref{KT_a} and \ref{KT_pt}.
\begin{figure}[h!]
\begin{center}
\subfigure[]{\includegraphics[width=0.45\textwidth,clip=]{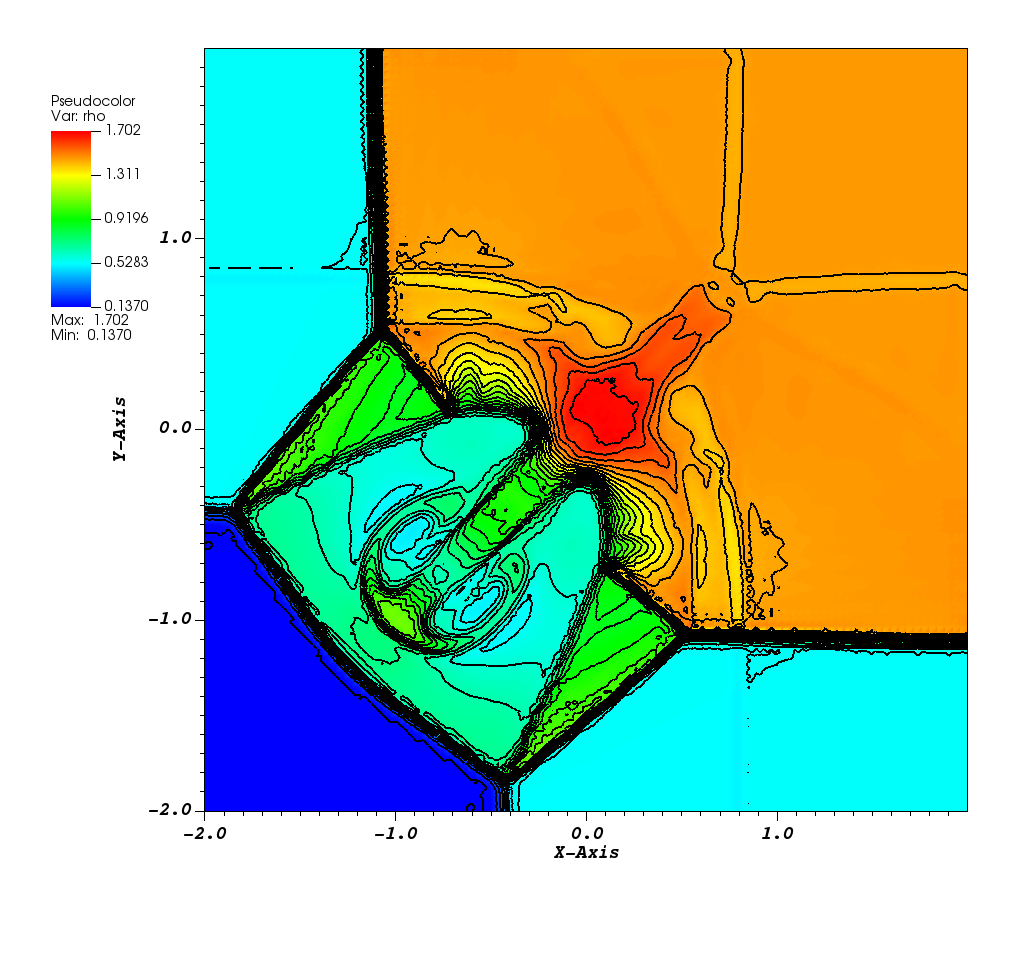}}
\subfigure[]{\includegraphics[width=0.45\textwidth,clip=]{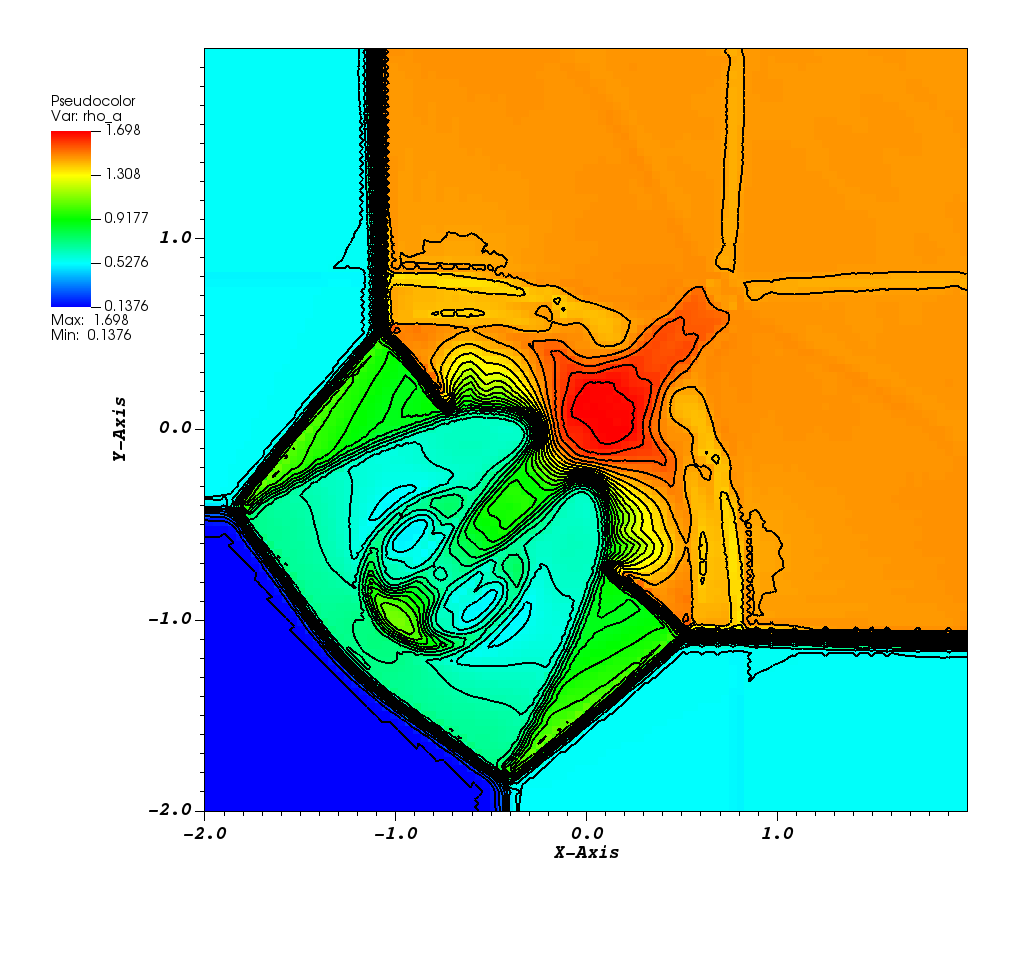}}
\end{center}
\caption{\label{KT_a} 30 equi-spaced isolines of the density as in Figure \ref{KT_old}, cfl=0.4 on a $100\times 100$ mesh. (a): averaged values, (b): point values.}
\end{figure}
\begin{figure}[h!]
\begin{center}
\subfigure[]{\includegraphics[width=0.45\textwidth,clip=]{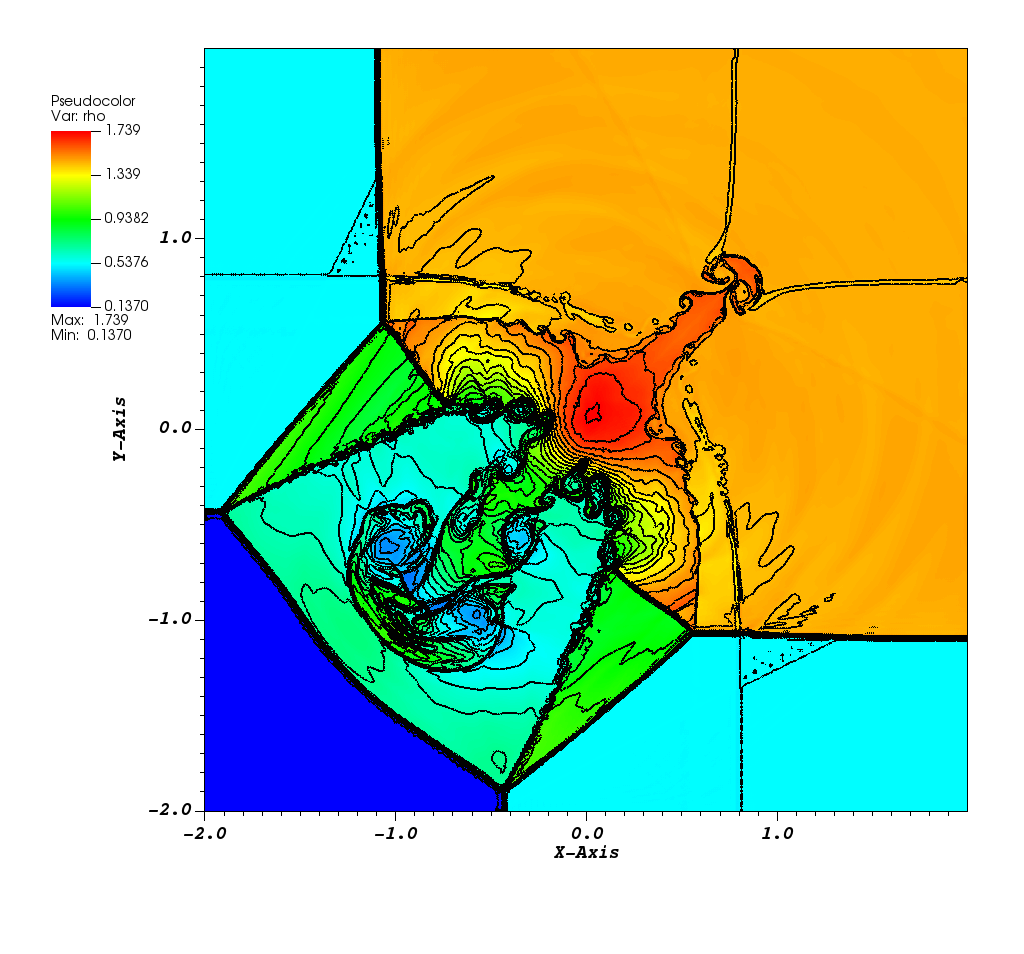}}
\subfigure[]{\includegraphics[width=0.45\textwidth,clip=]{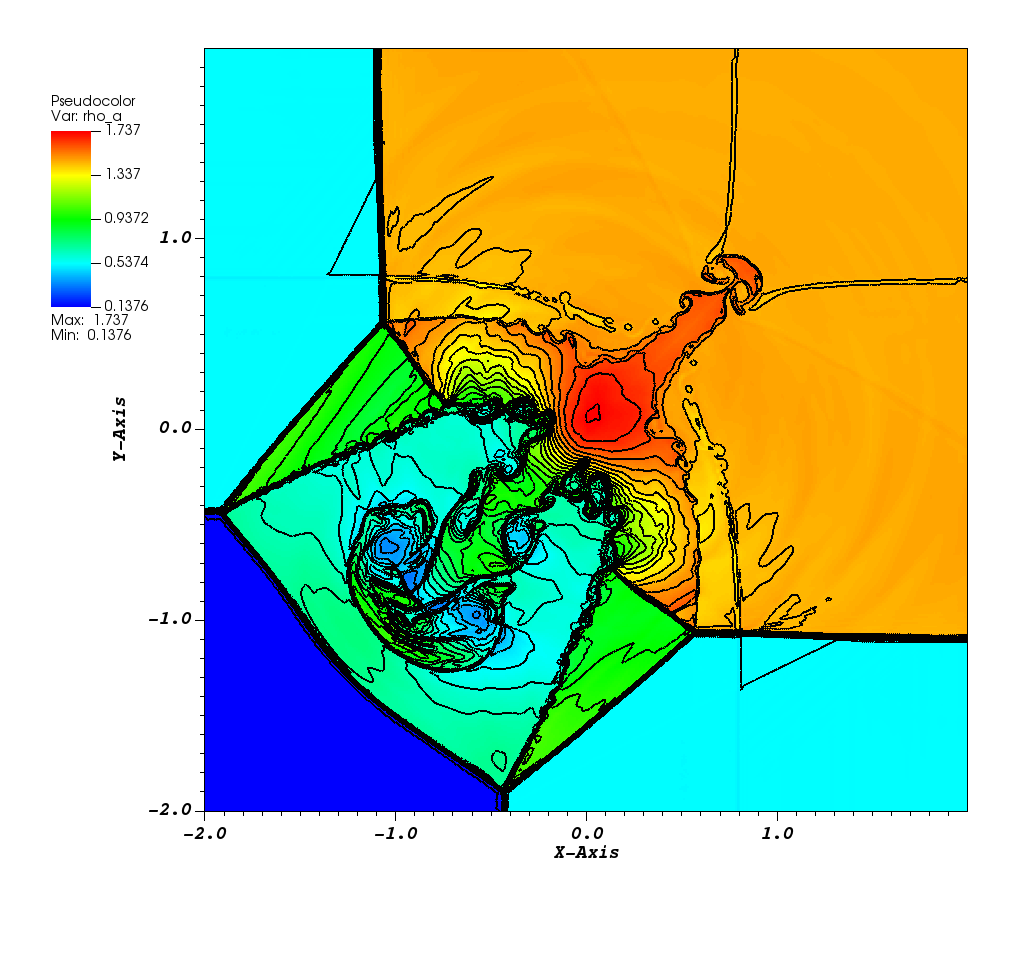}}
\end{center}
\caption{\label{KT_pt} 30 equi-spaced isolines of the density as in Figure \ref{KT_old}, cfl=0.4 on a $400\times 400$ mesh. (a): averaged values, (b): point values.}
\end{figure}
On Figure \ref{KT_old}, we have represented the solution obtained with the same $400\times400$ mesh, but with the algorithm of \cite{abgrall2024virtualfiniteelementhyperbolic} which uses a MOOD approach to stabilize the shocks.
\begin{figure}[h!]
\begin{center}
  \includegraphics[width=0.45\textwidth,clip=]{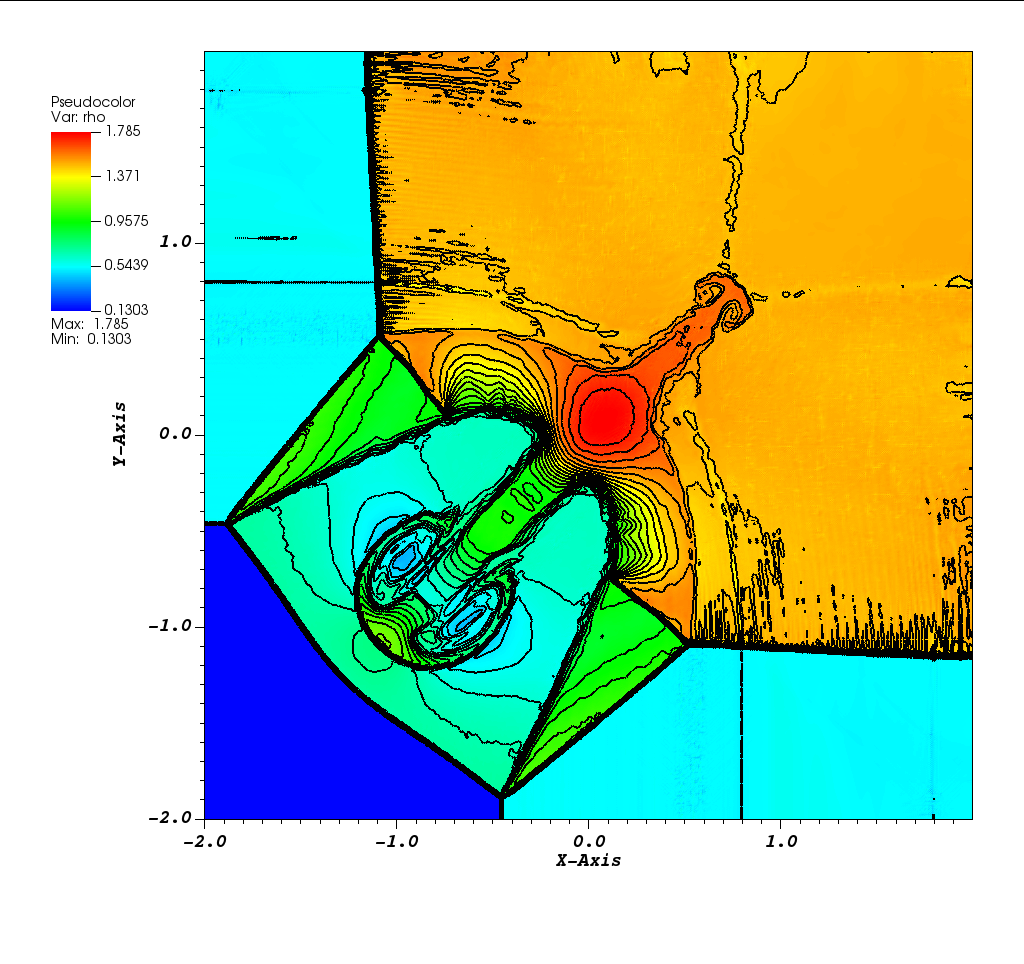}  
\end{center}
\caption{\label{KT_old}30 isolines $\rho\in [0.130,1.785]$, High order scheme \eqref{res:HO}-\eqref{flu:HO} with MOOD stabilisation, see \cite{abgrall2024virtualfiniteelementhyperbolic}.}
\end{figure}
We see that the new algorithm is way more accurate, and the solution is much cleaner. 
\section{Conclusion}
We have sketched a new method, inspired by the Roe's active flux method. The method works on general polygons, and is rigorously bound preserving. We have, because of space limitation, only presented results with shocks and on triangular meshes. Smooth solutions for scalar and the Euler equations are given in \cite{abgrall2024virtualfiniteelementhyperbolic}. We get the expected order of accuracy (i.e. 3).
%%%%%%%%%%%%%%%%%%%
\section*{Acknowledgement} We acknowledge many discussions with  Professor Kailiang Wu (SUSTECH, China).
\bibliographystyle{unsrt}
\bibliography{papier}

\begin{thebibliography}{10}

\bibitem{BPPampa}
R\'emi Abgrall, Yongle Liu, and Walter Boscheri.
\newblock Bound preserving virtual finite element and hyperbolic problems: the
  {PAMPA} algorithm.
\newblock {\em in preparation}, 2025.

\bibitem{abgr}
R{\'e}mi Abgrall.
\newblock The notion of conservation for residual distribution schemes (or
  fluctuation splitting schemes), with some applications.
\newblock {\em Commun. Appl. Math. Comput.}, 2(3):341--368, 2020.

\bibitem{AF1}
T.~A. Eyman and P.~L. Roe.
\newblock Active flux.
\newblock 49th AIAA Aerospace Science Meeting, 2011.

\bibitem{AF2}
T.~A. Eyman and P.~L. Roe.
\newblock Active flux for systems.
\newblock 20 th AIAA Computational Fluid Dynamics Conference, 2011.

\bibitem{AF3}
T.~A. Eyman.
\newblock {\em Active flux}.
\newblock PhD thesis, University of Michigan, 2013.

\bibitem{RoeAF}
P.L. Roe.
\newblock Is discontinuous reconstruction really a good idea?
\newblock {\em Journal of Scientific Computing}, 73:1094–1114, 2017.

\bibitem{VLeV}
B.~van Leer.
\newblock Towards the ultimate conservative difference scheme. {V}. {A}
  second-order sequel to {G}odunov's method.
\newblock {\em J. Comput. Phys.}, 32(1):101--136, 1979.

\bibitem{He}
Fanchen He and P.L. Roe.
\newblock A new treatment of conservation laws for a family of arbitrary-order
  fully-discrete numerical schemes based on the active flux method.
\newblock {\em Journal of Scientific Computing}, 2022.
\newblock submitted.

\bibitem{HKS}
C.~Helzel, D.~Kerkmann, and L.~Scandurra.
\newblock A new {ADER} method inspired by the active flux method.
\newblock {\em J. Sci. Comput.}, 80:35--61, 2019.

\bibitem{zbMATH07818681}
Erik Chudzik, Christiane Helzel, and M{\'a}ria
  Luk{\'a}{\v{c}}ov{\'a}-Medvi{\v{d}}ov{\'a}.
\newblock Active flux methods for hyperbolic systems using the method of
  bicharacteristics.
\newblock {\em J. Sci. Comput.}, 99(1):39, 2024.
\newblock Id/No 16.

\bibitem{zbMATH07698889}
Donna Calhoun, Erik Chudzik, and Christiane Helzel.
\newblock The {Cartesian} grid active flux method with adaptive mesh
  refinement.
\newblock {\em J. Sci. Comput.}, 94(3):31, 2023.
\newblock Id/No 54.

\bibitem{Abgrall_camc}
R.~Abgrall.
\newblock A combination of {R}esidual {D}istribution and the {A}ctive {F}lux
  formulations or a new class of schemes that can combine several writings of
  the same hyperbolic problem: application to the 1{D} {E}uler equation.
\newblock {\em Commun. Appl. Math. Comput.}, 5:370--402, 2023.

\bibitem{AB_HOAF}
R.~Abgrall and W.~Barsukow.
\newblock Extensions of {A}ctive {F}lux to arbitrary order of accuracy.
\newblock {\em ESAIM: Mathematical Modelling and Numerical Analysis},
  57:991--1027, 2023.

\bibitem{barsukow2025generalizedactivefluxmethod}
Wasilij Barsukow, Praveen Chandrashekar, Christian Klingenberg, and Lisa
  Lechner.
\newblock A generalized active flux method of arbitrarily high order in two
  dimensions, 2025.
\newblock arXiv 2502.05101.

\bibitem{zbMATH07695228}
Wasilij Barsukow.
\newblock Stationarity preservation properties of the active flux scheme on
  {Cartesian} grids.
\newblock {\em Commun. Appl. Math. Comput.}, 5(2):638--652, 2023.

\bibitem{abgrall2024virtualfiniteelementhyperbolic}
Rémi Abgrall, Yongle Liu, and Walter Boscheri.
\newblock Virtual finite element and hyperbolic problems: the {PAMPA}
  algorithm, 2024.
\newblock arXiv:2412.01341.

\bibitem{abgrall2024activefluxtriangularmeshes}
R\'emi Abgrall, Jianfang Lin, and Yongle Liu.
\newblock Active flux for triangular meshes for compressible flows problems.
\newblock {\em Beijing Journal of Pure and Applied Mathematics}, in press,
  2025.
\newblock arXiv 2312.11271.

\bibitem{Mood1}
S.~Clain, S.~Diot, and R.~Loub\`ere.
\newblock A high-order finite volume method for systems of conservation
  laws---{M}ulti-dimensional {O}ptimal {O}rder {D}etection ({MOOD}).
\newblock {\em J. Comput. Phys.}, 230(10):4028--4050, 2011.

\bibitem{wu2023geometric}
Kailiang Wu and Chi-Wang Shu.
\newblock Geometric quasilinearization framework for analysis and design of
  bound-preserving schemes.
\newblock {\em SIAM Review}, 65(4):1031--1073, 2023.

\bibitem{abgrall2024boundpreservingpointaveragemomentpolynomialinterpreted}
Rémi Abgrall, Miaosen Jiao, Yongle Liu, and Kailiang Wu.
\newblock Bound preserving point-average-moment polynomial-interpreted
  ({PAMPA}) scheme: one-dimensional case, 2024.

\bibitem{Kuzmin_MCL1}
D.~Kuzmin.
\newblock Monolithic convex limiting for continuous finite element
  discretizations of hyperbolic conservation laws.
\newblock {\em Comput. Method. Appl. M.}, 361:112804, 2020.

\bibitem{Zhang_MP}
X.~Zhang, Y.~Xia, and C.-W. Shu.
\newblock Maximum-principle-satisfying and positivity-preservinghigh order
  discontinuous galerkin schemes for conservation laws on triangular meshes.
\newblock {\em J. Sci. Comput.}, 50:29--62, 2012.

\bibitem{Guermond_IDP3}
J.-L. Guermond, B.~B.~Popov, and I.~Tomas.
\newblock Invariant domain preserving discretization independent schemes and
  convex limiting for hyperbolic systems.
\newblock {\em Comput. Method. Appl. M.}, 347:143--175, 2019.

\bibitem{duan2024AF}
J.~Duan, W.~Barsukow, and C.~Klingenberg.
\newblock Active flux methods for hyperbolic conservation laws--flux vector
  splitting and bound-preservation: One-dimensional case, 2024.
\newblock arXiv 2405.02447.

\bibitem{BPPAMPA1D}
R.~Abgrall, M.~Jiao, Y.~Liu, and K.~Wu.
\newblock Bound preserving {P}oint-{A}verage-{M}oment
  {P}olynomi{A}l-interpreted ({PAMPA}) scheme: one-dimensional case, 2024.
\newblock arXiv:2410.14292.

\bibitem{hitch}
L.~Beir{\~a}o~da Veiga, Franco Brezzi, L.~D. Marini, and A.~Russo.
\newblock The {Hitchhiker}'s guide to the virtual element method.
\newblock {\em Math. Models Methods Appl. Sci.}, 24(8):1541--1573, 2014.

\bibitem{GuermondPopovFast}
Jean-Luc Guermond and Bojan Popov.
\newblock Fast estimation from above of the maximum wave speed in the {Riemann}
  problem for the {Euler} equations.
\newblock {\em J. Comput. Phys.}, 321:908--926, 2016.

\bibitem{GauthierLBM}
G.~Wissocq, Y.~Liu, and R.~Abgrall.
\newblock A positive- and bound-preserving vectorial lattice boltzmann method
  in two dimensions.
\newblock {\em SIAM J. Sci. Computing}, 2025.
\newblock accepted.

\end{thebibliography}
\end{document}